\documentclass[12pt, twoside]{article}
\usepackage{pcam}

\def\d{{\rm d}}
\def\E{{ e}}
\def\i{{ i}}

\makeatletter
\newcommand{\mathleft}{\@fleqntrue\@mathmargin0pt}
\newcommand{\mathcenter}{\@fleqnfalse}
\makeatother

\newtheorem{theorem}{Theorem}

\begin{document}
\English

\title
	[Optimal quadrature formula] 
    {Optimal quadrature formulas for computing of Fourier integrals in a Hilbert space} 
\author
	[Hayotov~A.R.] 
	{Hayotov~A.R.,Babaev~S.S. } 
    [$^{1,3}$Hayotov~A.R., $^{1,2}$Babaev~S.S. ] 
\email
    { hayotov@mail.ru; b\_samandar@mail.ru}
\thanks
 {The work has been done while Samandar S. Babaev was visiting Department of Mathematical
Sciences at KAIST, Daejeon, Republic of Korea, as a fellow of the 'El yurt umidi'(EYUF) Foundation of Uzbekistan  for advanced training course.
}
\organization
    {$^1$V.I.Romanovskiy Institute of Mathematics, Uzbekistan Academy of Sciences, 81, M.Ulugbek str., Tashkent 100170, Uzbekistan;

    $^2$Bukhara State University, 11, M.Ikbol str., Bukhara 200114, Uzbekistan;

     $^3$National University of Uzbekistan named after Mirzo Ulugbek, 4, University str., Tashkent 100174, Uzbekistan
     }
\abstract
{In the present paper the optimal quadrature formulas in the sense of Sard are constructed for numerical integration of the integral $\int_a^b\E^{2\pi\i\omega x}\varphi(x)\d x$ with $\omega\in \mathbb{R}$ in the Hilbert space $W_2^{(2,1)}[a,b]$ of complex-valued functions. Furthermore, the explicit expressions for coefficients of the constructed optimal quadrature formulas are obtained.
At the end of the paper some numerical results are presented.}


\titleRus
    [Оптимальная квадратурная формула] 
    {Оптимальные квадратурные формулы для вычисления интегралов Фурье в гильбертовом пространстве} 
\authorRus
    [Хаётов ~ А.Р.] 
    {Хаётов ~ А.Р., Бабаев ~ С.С.} 
    [ $ ^ {1,3} $Хаётов ~ А.Р., $ ^ {1,2} $ Бабаев ~ С.С.] 
\thanksRus
    {Работа была проделана в то время, когда Самандар С. Бабаев посещал Отдел математических наук в KAIST, Тэджон, Республика Корея, в качестве стипендиата Фонда «Эл-юрт умиди».}
\organizationRus
    {$ ^ 1 $Институт математики им. В.И.Романовского Академии наук Узбекистана,  ул. М.Улугбека, 81, Ташкент 100170, Узбекистан;

$ ^ 2 $Бухарский государственный университет, ул. М.Икбола, 11, Бухара, 200114, Узбекистан;

$ ^ 3 $Национальный университет Узбекистана имени Мирзо Улугбека, ул. Университетская, 4, Ташкент 100174, Узбекистан
}
\abstractRus
    { В настоящей статье построены оптимальные квадратурные формулы в смысле Сарда для приближенного вычисления интеграла  $\int_a^b\E^{2\pi\i\omega x}\varphi(x)\d x$ с  $\omega\in \mathbb{R}$
в гильбертовом пространстве $W_2^{(2,1)}[a,b]$ комплекснозначных функций. Кроме того, получены явные выражения для коэффициентов построенных
оптимальных формул. В конце статье приведены некоторые численные результаты.}

\UDC{519.644}
\receivedRus{31.03.2020}
\receivedEng{March 31, 2020}

\maketitle

\setcounter{section}{0}
\section{Introduction and statement of the Problem}

We wish to find efficient numerical approximations for the Fourier integrals of the form
\begin{equation}\label{eq1.1}
I(\varphi)=\int\limits_a^b\E^{2\pi \i \omega x}\varphi(x)\d x
\end{equation}
with $\omega\in \mathbb{R}$. This type of integrals are called \emph{highly oscillating integrals}.
In most cases it is impossible to get the exact values of such integrals.
Thus, they can be approximately calculated using the formulas of numerical integration.
However, standard methods of numerical integration cannot be successfully applied for that.
Therefore special effective methods should be developed for approximation of highly oscillating integrals.
One of the first numerical integration formula for the integral (\ref{eq1.1}) was obtained by Filon \cite{Filon28} in 1928 using a quadratic spline.
Since then, for integrals of different types of highly oscillating functions many special effective methods have been developed, such as
Filon-type method, Clenshaw-Curtis-Filon type method, Levin type methods, modified Clenshaw-Curtis method, generalized quad\-rature rule, and Gauss-Laguerre quadrature (see, for example, \cite{AvdMal89,IBab,BabVitPrag69,BakhVas68,IserNor05,Mil98,NovUllWoz15,ZhangNovak19}, for more review see, for instance, \cite{DeanoHuyIser18,MilStan14,Olver08} and references therein).

Recently, in \cite{BolHayShad17,BolHayMilShad17}, based on Sobolev's method, the problem of construction of optimal quadrature formulas in the sense of Sard for numerical calculation of the integral (\ref{eq1.1}) with integer $\omega$ was studied in Hilbert spaces $L_2^{(m)}$ and $W_2^{(m,m-1)}$, respectively.
The optimal quadrature formulas are constructed in \cite{Samandar1,HayJeonLee19} for approximation Fourier integrals in the spaces $W_2^{(1,0)}[a,b]$ and $L_2^{(1)}$, respectively.

Here, we consider the Hilbert space $W_2^{(2,1)}[a,b]$ of non-periodic, complex-valued functions defined on the interval $[a,b]$ which posses an absolutely continuous first derivative on $[a,b]$, and whose second order derivative is in $L_2$ and in this space the inner product is defined by the equality
\begin{equation*}
\langle\varphi,\psi\rangle=\int\limits_a^b(\varphi''(x)+\varphi'(x))(\bar\psi''(x)+\bar\psi'(x))\d x,
\end{equation*}
where $\bar\psi$ is the conjugate function to the function $\psi$ and the norm of the function $\varphi$ is correspondingly defined by the formula
$$
\|\varphi\|_{W_2^{(2,1)}[a,b]}=\langle\varphi,\varphi\rangle^{\frac{1}{2}}
$$
and $\int\limits_a^b(\varphi''(x)+\varphi'(x))(\bar\varphi''(x)+\bar\varphi'(x))\d x<\infty$.

We consider the following quadrature formula
\begin{equation}
\int\limits_a^b {\E^{2\pi i\omega x}\varphi (x)\d x \cong } \sum\limits_{\beta  = 0}^N
{C_\beta  } \varphi (h\beta  )\label{eq1.3}
\end{equation}
with the error functional
\begin{equation}
\ell (x) =\E^{2\pi i\omega x}\varepsilon _{[a,b]} (x) - \sum\limits_{\beta  = 0}^N
{C_\beta  } \delta (x - h\beta  ),\label{eq1.4}
\end{equation}
where $C_\beta$ are the coefficients of formula (\ref{eq1.3}), $h=\frac{b-a}{N}$,
$N\in \mathbb{N}$, $i^2=-1$, $\omega\in \mathbb{R}$,  $\varepsilon _{[a,b]} (x)$ is the indicator of the interval $[a,b]$ and $\delta(x)$ is the Dirac delta-function.
We mention that the coefficients $C_{\beta}$ depend on $\omega$ and $h$.

The difference
\begin{equation*}
(\ell,\varphi) = \int\limits_a^b {\E^{2\pi \omega ix}\varphi (x)dx -
\sum\limits_{\beta  = 0}^N {C_\beta  \varphi (x_\beta  )} }  =
\int\limits_{ - \infty }^\infty  {\ell (x)\varphi
(x)\d x}
\end{equation*}
is called \emph{the error} of the quadrature formula
(\ref{eq1.3}).

The error functional (\ref{eq1.4}) is a linear functional
in $W_2^{(2,1)*}[a,b]$, where $W_2^{(2,1)*}[a,b]$ is the
conjugate space to the space $W_2^{(2,1)}[a,b]$. Since the functional (\ref{eq1.4}) is defined on the space $W_2^{(2,1)}$ the conditions
\begin{equation}\label{eq1.6}
(\ell,1)=0,
\end{equation}
\begin{equation}\label{eq1.7}
(\ell,\E^{-x}) = 0
\end{equation}
should be fulfilled. The equations (\ref{eq1.6}) and (\ref{eq1.7}) mean the exactness of the quadrature formula (\ref{eq1.3}) for any linear combination of functions $1$ and $\E^{-x}$.

Here we consider the problem of construction of optimal quadrature formulas of the (\ref{eq1.3}) in the sense of Sard in the space $W_2^{(2,1)}[a,b]$.
Sard's optimization problem of numerical integration formulas of the form (\ref{eq1.3}) in the space $W_2^{(2,1)}[a,b]$
is the problem of finding the minimum of the norm of the error functional $\ell$ by coefficients $C_{\beta}$, i.e.,
to find coefficients $C_{\beta}$ satisfying the equality
\begin{equation}\label{eq1.8}
\|\mathring{\ell}\|_{W_2^{(2,1)*}[a,b]}=\inf\limits_{C_\beta}\|\ell\|_{W_2^{(2,1)*}[a,b]}.
\end{equation}
The coefficients satisfying the last equality are called \emph{optimal one} and they are denoted as $\mathring{C}_{\beta}$.
The  quadrature formula with coefficients $\mathring{C}_{\beta}$ is called \emph{the optimal quadrature formula in the sense of Sard}, and $\mathring{\ell}$ is the error functional corresponding to the optimal quadrature formula.
The solution of Sard's problem gives the sharp upper bound for the error $(\ell,\varphi)$ of functions $\varphi$ from the space $W_2^{(2,1)}[a,b]$ as follows
$$
|(\ell,\varphi)|\leq \|\varphi\|_{W_2^{(2,1)}[a,b]}
\|\mathring{\ell}\|_{W_2^{(2,1)*}[a,b]}.
$$

This problem, for the quadrature formulas of the form (\ref{eq1.3}) with $\omega=0$ and $\omega \in \mathbb{Z}$ with $\omega\neq 0$, was solved in the works \cite{ShadHay14,BolHayMilShad17} in the space $W_2^{(m,m-1)}$ for any integer $m$.

Further, we solve Sard's problem on construction of optimal quadrature formulas of the form (\ref{eq1.3}) for $\omega\in \mathbb{R} $ with $\omega\neq 0$, first for the interval $[0,1]$ and then using a linear transformation for the interval $[a,b]$.
For this we use the following auxiliary results.

We need the concept of discrete argument functions and operations on them (see, for instance, \cite{Sobolev06,Sobolev74}).

In finding the analytic formulas for coefficients of optimal formulas in the space $W_2^{(2,1)}$ by Sobolev method the discrete analogue $D_2(h\beta)$ of the operator ${{\d^{4}}\over {\d x^{4}}}-{{\d^{2}}\over {\d x^{2}}}$ \cite{ShadHay04a,ShadHay04b} plays the main role.
This discrete analogue  satisfies the equality
\begin{equation}\label{eq1.9}
D_2(h\beta)*G_2(h\beta) = \delta_{\d}(h\beta),
\end{equation}
where $G_2(h\beta)$ is the discrete argument function for the function
\begin{equation}\label{eq1.10}
G_2(x) = {{{\rm{sgn}}x} \over2}\left(\frac{e^x-e^{-x}}{2}-x\right)
\end{equation}
and $\delta_{\d}(h\beta)$ is equal to 0 for $\beta\neq 0$ and 1 for $\beta=0$.

\medskip

\begin{theorem}\label{Thm1}  The discrete analogue to the
differential operator ${{\d^{4}} \over {\d x^{4} }}-{{\d^{2}}
\over {\d x^{2} }}$ satisfying the equation (\ref{eq1.9}) has the form
\begin{equation}
 D_2(h\beta)=\frac{1}{p}\left\{
\begin{array}{ll}
A \lambda_1^{|\beta|-1},& |\beta|\geq 2,\\
-2\E^h+A,& |\beta|=1,\\
2C+\frac{A}{\lambda_1},& \beta=0,
\end{array}
\right.\label{eq1.11}
\end{equation}
where
\begin{eqnarray*}
C&=&1+2\E^h+\E^{2h}-\frac{\E^h\cdot (\lambda_1^2+1)}{\lambda_1},
\\
A&=&\frac{2(\lambda_1-1)[\lambda_1(\E^{2h}+1)-\E^h(\lambda_1^2+1)]}{\lambda_1+1} ,
\\
\mathcal{P}_{2}(\lambda)&=&p\lambda^2-2(1-\E^{2h}+h(\E^{2h}+1))\lambda+p,\\
p&=&1+2h\E^h-\E^{2h},\\
\end{eqnarray*}
$\lambda_1$ is the root of the polynomial $\mathcal{P}_{2}(\lambda)$ for which $|\lambda_1|<1.$
\end{theorem}

In addition, some properties of $D_2(h\beta)$ were studied in the works \cite{ShadHay04a,ShadHay04b}. Here we give the
following.

\begin{theorem}\label{Thm2} { The discrete analogue $D_2(h\beta)$ to
the differential operator ${{\d^{4}} \over {\d x^{4}}}-{{\d^{2}} \over {\d x^{2}}}$ satisfies the following
equalities

1) $D_2(h\beta)*\E^{h\beta}=0,$

2) $D_2(h\beta)*\E^{-h\beta}=0,$

3) $D_2(h\beta)*(h\beta)^n=0,$\  \ for\  \ $n=0,1$,

4) $D_2(h\beta)*G_2(h\beta)=\delta_{\mathrm{d}}(h\beta),$\\
here $G_2(h\beta)$ is the function of discrete argument
corresponding to the function $G_2(x)$, defined by equality
(\ref{eq1.10}) and $\delta_{\mathrm{d}}(h\beta)$ is the discrete delta
function.}
\end{theorem}

Now we give the results of the work \cite{ ShadHay14} on the optimal quadrature formulas
of the form (\ref{eq1.3}) in the sense of Sard for the case $\omega=0$.

 Coefficients for the optimal quadrature formula
\begin{equation}\nonumber
\int\limits_0^1\varphi(x)\d x\cong \sum\limits_{\beta=0}^NC_{\beta}\varphi(h\beta)
\end{equation}
in the space $W_2^{(2,1)}[0,1]$ have the form
\begin{eqnarray} \label{eq1.13}
\mathring{C}_0&=&1-\frac{h}{e^h-1}-K(\lambda_1-\lambda_1^N),\nonumber\\
\mathring{C}_{\beta}&=&h+K\left((e^h-\lambda_1)\lambda_1^{\beta}+(1-\lambda_1e^h)\lambda_1^{N-\beta}\right),\
\beta=1,...,N-1,\\
\mathring{C}_N&=&-1+\frac{e^hh}{e^h-1}-K(\lambda_1-\lambda_1^N)e^h,\nonumber
\end{eqnarray}
where
\begin{eqnarray}
K&=&\frac{(2e^h-2-he^h-h)(\lambda_1-1)}{2(e^h-1)^2(\lambda_1+\lambda_1^{N+1})},\nonumber\\
\lambda_1&=&\frac{h(e^{2h}+1)-e^{2h}+1-(e^h-1)\sqrt{h^2(e^{h}+1)^2+2h(1-e^h)}}
{1-e^{2h}+2he^h}\nonumber \\
&=&\sqrt{3}-2+O(h^2),\label{eq1.12}
\end{eqnarray}
$ |\lambda_1|<1$, $h=1/N$, $N=2,3,...$.

\section{Construction of optimal quadrature formulas for the interval $[0,1]$}

Here we obtain optimal quadrature formulas of the form (\ref{eq1.3}) for the interval $[0,1]$ when $\omega\in \mathbb{R}$ and $\omega\neq 0$. In the space $W_2^{(2,1)}[0,1]$ for the coefficients of the optimal quadrature formula in the sense of Sard of the form
\begin{equation}\nonumber
\int\limits_0^1 {\E^{2\pi \i\omega x}\varphi (x)\d x \cong } \sum\limits_{\beta  = 0}^N
{C_\beta  } \varphi (h\beta  )
\end{equation}
for $\omega\in \mathbb{R}$ with $\omega\neq 0$, we get the following system of linear equations

\begin{eqnarray}
&&\sum\limits_{\gamma=0}^N C_\gamma G_2(h\beta   - h\gamma) +p_0+d\E^{-h\beta} =f_2(h\beta), \ \beta=0,1,...,N, \label{eq2.1}\\
&&\sum\limits_{\gamma=0}^N C_\gamma=\frac{\E^{2\pi i \omega}-1}{2\pi i \omega},\label{eq2.2}\\
&&\sum\limits_{\gamma=0}^N C_\gamma\E^{-h\gamma}=\frac{\E^{2\pi i \omega-1}-1}{2\pi i \omega-1},\label{eq2.3}
\end{eqnarray}
where $p_0$ is a complex number,
\begin{eqnarray}\nonumber
f_2(h\beta)&=&\frac{\E^{-h\beta}}{4}\cdot \frac{\E^{2\pi i\omega+1}-2\E^{(2\pi i\omega+1)h\beta}+1}{2\pi i\omega+1}-\frac{\E^{h\beta}}{4}\cdot \frac{\E^{2\pi i\omega-1}-2\E^{(2\pi i\omega-1)h\beta}+1}{2\pi i\omega-1}\\
&&+\frac{\E^{2\pi i\omega}-2\E^{2\pi i\omega h\beta}+1}{2(2\pi i\omega)^2}+\frac{h\beta\E^{2\pi i\omega}+h\beta-\E^{2\pi i\omega}}{4\pi i\omega},\label{eq2.4}
\end{eqnarray}
$G_2(x)$ is defined by (\ref{eq1.10}), $h=1/N$, $N$ is a natural number.

In the system (\ref{eq2.1})-(\ref{eq2.3}) unknowns are the coefficients $C_\beta,$ $\beta=0,1,...,N$, $p_0$ and $d$. We point out that the system (\ref{eq2.1})-(\ref{eq2.3}) has a unique solution.

This solution satisfies conditions (\ref{eq1.6}), (\ref{eq1.7}) and equality (\ref{eq1.8}). It should be noted that
the existence and uniqueness of the solution for such type of systems were studied, for example,
in \cite{Sobolev74}.

We are interested in finding explicit expressions for the optimal coefficients $\mathring{C}_{\beta}$, $\beta=0,1,...,N$ and unknowns  $p_0$ and $d$ satisfying the system (\ref{eq2.1})-(\ref{eq2.3}). The system (\ref{eq2.1})-(\ref{eq2.3}) is solved similarly as the system (34)-(35) of \cite{ShadHay14}  using the discrete analogue $D_2(h\beta)$.

We formulate the results of this section as the following two theorems.

\begin{theorem}\label{Thm6}
{For real $\omega$ with $\omega h\not\in \mathbb{Z}$, the coefficients of optimal quadrature
formulas of the form (\ref{eq1.3})  in the space $W_2^{(2,1)}[0,1]$  are expressed by formulas}
\begin{eqnarray}
\mathring{C}_0&=&\frac{K_{\omega,2}\ \E^{2\pi i\omega h}}{\E^{2\pi i\omega h}-1}
-\frac{1}{2\pi i\omega}+a_1\frac{\lambda_1}{\lambda_1-1}+b_1\frac{\lambda_1^{N}}{1-\lambda_1}, \label{eq2.5} \\
\mathring{C}_\beta&=&\E^{2\pi \i\omega h\beta}K_{\omega,2}+a_1\lambda_1^\beta+b_1\lambda_1^{N-\beta},
\qquad\qquad \beta=1,2,...,N-1, \label{eq2.6}\\
\mathring{C}_N&=&K_{\omega,2}\ \Bigg(\frac{\E^{2\pi i\omega}\E^{h}}{\E^{h}-\E^{2\pi i\omega h}}+\frac{\E^{2\pi i\omega h+1}(1-\E^h)}{(\E^{2\pi i\omega h}-1)(\E^h-\E^{2\pi i\omega h})} \Bigg) \nonumber\\
&&+\frac{\E^{2\pi i\omega}}{2\pi i\omega-1}+a_1\frac{\lambda_1^{N}\E^{h}}{\E^{h}-\lambda_1} +b_1\frac{\lambda_1\E^{h}}{\lambda_1\E^{h}-1},  \label{eq2.7}
\end{eqnarray}

{where }\\
\begin{eqnarray}\nonumber
&&a_1=\frac{(\lambda_1-1)(\lambda_1-\E^h)}{\lambda_1(1-\lambda_1^{2N})}\Bigg(\frac{1-\lambda_1^N\E^{2\pi i\omega}}{2\pi i\omega(2\pi i\omega-1)(\E^h-1)}+\frac{K_{\omega,2}\lambda_1^N(\E-\E^{2\pi i\omega})\E^{2\pi i\omega h}}{(\E^{2\pi i\omega h}-1)(\E^h-\E^{2\pi i\omega h})}\Bigg),\\
&&b_1= \frac{(\lambda_1-1)(1-\lambda_1\E^h)}{\lambda_1(1-\lambda_1^{2N})}\Bigg(\frac{\lambda_1^N-\E^{2\pi i\omega}}{2\pi i\omega(2\pi i\omega-1)(\E^h-1)}+\frac{K_{\omega,2}(\E-\E^{2\pi i\omega})\E^{2\pi i\omega h}}{(\E^{2\pi i\omega h}-1)(\E^h-\E^{2\pi i\omega h})}\Bigg)\label{eq2.8}\\ \nonumber
\end{eqnarray}
and
\begin{equation} \label{eq2.9}
K_{\omega,2}=\frac{ (\E^h-\E^{2\pi i\omega h})(1-\E^{2\pi i\omega h+h})(1-\E^{2\pi i\omega h})^2}{2\pi^2\omega^2(4\pi^2\omega^2+1)\E^{2\pi i\omega h}\left((1-\E^{2h})(1-\E^{2\pi i\omega h})^2+2h(\E^h-\E^{2\pi i\omega h})(1-\E^{2\pi i\omega h+h})\right)},
\end{equation}
$\lambda_1$ is defined by (\ref{eq1.12}), $h=1/N$, $N$ is a natural number.
\end{theorem}

\begin{theorem}\label{Thm7} {The coefficients of optimal quadrature formulas of the form (\ref{eq1.3}) with the error
functional (\ref{eq1.4})  when $\omega h\in \mathbb{Z}$ and $\omega\neq 0$ in the
space $W_2^{(2,1)}[0,1]$ are expressed by formulas
\begin{eqnarray*}
\mathring{C}_0&=&\frac{2\pi i\omega (1-\E^h)-1}{2\pi i\omega (1-2\pi i\omega)(1-\E^h)}+a_1\frac{\lambda_1^2}
{(1-\lambda_1)(\E^h-\lambda_1)}+b_1\frac{\lambda_1^N}
{(1-\lambda_1)(1-\lambda_1\E^h)},\\
\mathring{C}_\beta&=&a_1\lambda_1^\beta+b_1\lambda_1^{N-\beta},
\qquad\qquad\qquad \beta=1,2,...,N-1,\\
\mathring{C}_N&=&\frac{2\pi i\omega(\E^h-1)-\E^h}{2\pi i\omega(1-2\pi i\omega)(1-\E^h)}+ a_1\frac{\E^h\lambda_1^N}
{(1-\lambda_1)(\E^h-\lambda_1)}+b_1\frac{\E^h\lambda_1^2}
{(1-\lambda_1)(1-\lambda_1\E^h)},
\end{eqnarray*}
where
$$
a_1= \frac{(\E^h-\lambda_1)(1-\lambda_1)}{2\pi i \omega (2\pi i \omega -1)\lambda_1(\E^h-1)(\lambda_1^N+1)},
$$
$$
b_1=\frac{(1-\E^h\lambda_1)(1-\lambda_1)}{2\pi i \omega (2\pi i \omega -1)\lambda_1(\E^h-1)(\lambda_1^N+1)},
$$
$\lambda_1$ is defined by (\ref{eq1.12}), $h=1/N$, $N$ is a natural number.
}
\end{theorem}

We note that Theorem \ref{Thm6} is generalization of Theorem 6 in \cite{BolHayMilShad17} for real $\omega$ with $\omega h\not\in \mathbb{Z}$ while Theorem \ref{Thm7} for $\omega h\in \mathbb{Z}$ with $\omega\neq 0$ is the same with Theorem 7 of the work \cite{BolHayMilShad17}.
Therefore, it is sufficient to give a  proof of Theorem \ref{Thm6}.

\medskip

\emph{The proof of Theorem \ref{Thm6}}. First, as in Theorem 5 of the work \cite{ShadHay14}, using the discrete function $D_2(h\beta)$, for
optimal coefficients $\mathring{C}_{\beta},$ $\beta=1,2,...,N-1$, when $\omega$ is real and $\omega h\not\in \mathbb{Z}$ we get the following formula
\begin{equation}\label{eq2.10}
\mathring{C}_\beta=K_{\omega,2}\E^{2\pi\i \omega h\beta}+a_1\lambda_1^{\beta}+b_1\lambda_1^{N-\beta} ,\ \
\beta=1,2,...,N-1,
\end{equation}
where $a_1$, $b_1$, $K_{\omega,2}$ are unknowns,
and $\lambda_1$ is the root of the polynomial $\mathcal{P}_{2}(\lambda)$, $|\lambda_1|<1$.
Next, it is sufficient to find $a_1,\ b_1$, $K_{\omega,2}$,
 $\mathring{C}_0$, $\mathring{C}_N$ and unknown coefficients $p_0$ and $d$.

Now putting the form (\ref{eq2.10}) of optimal coefficients $\mathring{C}_{\beta}$, $\beta=1,2,...,N-1$ into (\ref{eq2.1}), (\ref{eq2.2}), (\ref{eq2.3}) and (\ref{eq2.4}) after some simplifications,
we get the following identity with respect to $(h\beta)$:
$$
\begin{array}{l}
\displaystyle \frac{\E^{h\beta}}{2}\left(C_0+\frac{K_{\omega,2}\E^{2\pi i\omega h}}{\E^h-\E^{2\pi i\omega h}}+a_1\frac{\lambda_1}{\E^h-\lambda_1}+b_1\frac{\lambda_1^N}{\E^h\lambda_1-1}-\frac{\E^{2\pi i\omega-1}-1}{2(2\pi i\omega-1)}\right)\\
\displaystyle  -\frac{\E^{-h\beta}}{2}\left(C_0+\frac{K_{\omega,2}\E^{h+2\pi i\omega h}}{1-\E^{h+2\pi i\omega h}}+a_1\frac{\lambda_1\E^h}{1-\E^h\lambda_1}+b_1\frac{\lambda_1^N\E^h}{\lambda_1-\E^h}-\frac{1}{2}\sum\limits_{\gamma=0}^N C_\gamma \E^{h\gamma}-2d\right)\\
\displaystyle  -\frac{\E^{2\pi i\omega h\beta}}{2}\left(\frac{K_{\omega,2}\E^{h}}{\E^h-\E^{2\pi i\omega h}}-\frac{K_{\omega,2}}{1-\E^{h+2\pi i\omega h}}+\frac{2hK_{\omega,2}\E^{2\pi i\omega h}}{(\E^{2\pi i\omega h}-1)^2}\right)\\
\displaystyle +(h\beta)\left(-C_0+\frac{K_{\omega,2}\E^{2\pi i\omega h}}{\E^{2\pi i\omega h}-1}+a_1\frac{\lambda_1}{\lambda_1-1}+b_1\frac{\lambda_1^N}{1-\lambda_1}+\frac{\E^{2\pi i\omega}-1}{4\pi i\omega}\right)
\end{array}
$$
$$
\begin{array}{l}
\displaystyle +\frac{hK_{\omega,2}\E^{2\pi i\omega h}}{(\E^{2\pi i\omega h}-1)^2}+ha_1\frac{\lambda_1}{(\lambda_1-1)^2}+hb_1\frac{\lambda_1^{N+1}}{(\lambda_1-1)^2}-\frac{h}{2}\sum\limits_{\gamma=0}^NC_\gamma \gamma+p_0\\
\displaystyle \qquad\qquad \qquad =-\frac{\E^{h\beta}}{2}\frac{\E^{2\pi i\omega-1}+1}{2(2\pi i\omega-1)}+\frac{\E^{-h\beta}}{2}\frac{\E^{2\pi i\omega+1}+1}{2(2\pi i\omega+1)}
-\frac{\E^{2\pi i\omega h\beta}}{2}\Bigg(\frac{1}{2\pi i\omega+1}\\
\displaystyle \qquad\qquad\qquad\quad-\frac{1}{2\pi i\omega-1}+\frac{2}{(2\pi i\omega)^2}\Bigg)+(h\beta)\frac{\E^{2\pi i\omega}+1}{4\pi i\omega}+\frac{\E^{2\pi i\omega}+1}{2(2\pi i\omega)^2}-\frac{\E^{2\pi i\omega}}{4\pi i\omega}
\end{array}
$$
Hence, the formula for $K_{\omega,2}$ is obtained by equating the coefficients of the term  $\E^{2\pi\i\omega h\beta}$.
Then we get the following system of two linear equations for $a_1$ and $b_1$. The first equality is found from equating the coefficients of the term $\E^{h\beta}$ and the second equality is found  from   equation (\ref{eq2.3})
\begin{equation}\label{eq2.11}
\begin{array}{ll}
 a_1 \frac{\lambda_1-\lambda_1^{N+1}}{(\lambda_1-1)(\lambda_1-\E^h)}+ b_1 \frac{\lambda_1-\lambda_1^{N+1}}{(\lambda_1-1)(1-\lambda_1\E^h)}=\frac{1-\E^{2 \pi i \omega}}{2 \pi i \omega(2 \pi i \omega-1)(\E^h-1)}+K_{\omega,2}\frac{(1-\E^{2\pi i \omega})\E^{2\pi i \omega h}}{(\E^{2\pi i \omega h}-1)(\E^h-\E^{2\pi i \omega h})},\\[2mm]
a_1 \frac{\lambda_1 \E-\lambda_1^{N+1}}{(\lambda_1-1)(\lambda_1-\E^h)}+ b_1 \frac{\lambda_1^{N+1} \E-\lambda_1}{(\lambda_1-1)(\lambda_1\E^h-1)}=\frac{\E-\E^{2 \pi i \omega}}{2 \pi i \omega(2 \pi i \omega-1)(\E^h-1)}+K_{\omega,2}\frac{(\E-\E^{2\pi i \omega})\E^{2\pi i \omega h}}{(\E^{2\pi i \omega h}-1)(\E^h-\E^{2\pi i \omega h})}.
\end{array}
\end{equation}
Also we find $C_0$ by equating the coefficients of the term $(h\beta)$, coefficient $C_N$ is easily found from  (\ref{eq2.2}).
 Finally, $a_1$ and $b_1$ are found from (\ref{eq2.11}).  Theorem  \ref{Thm6} is proved. \hfill $\Box$

\section{Optimal quadrature formulas for the interval $[a,b]$}

In this section we obtain the optimal quadrature formulas in the interval $[a,b]$ in the space $W_2^{(2,1)}$
by a linear transformation from the results of the previous section.

We consider construction of the optimal quadrature formulas of the form
\begin{equation}\label{eq2.12}
\int\limits_a^b\E^{2\pi i \omega x}\varphi(x)\ \d x\cong \sum\limits_{\beta=0}^NC_{\beta,\omega}[a,b]\varphi(x_\beta)
\end{equation}
in the Hilbert space $W_2^{(2,1)}[a,b]$. Here $C_{\beta,\omega}[a,b]$ are coefficients,
$x_\beta=h\beta+a$ $(\in [a,b])$ are nodes of the formula (\ref{eq2.12}),
$\omega\in \mathbb{R}$, $i^2=-1$ and $h=\frac{b-a}{N}$.

Now, by a linear transformation $x=(b-a)y+a$, where $0\leq y\leq 1$, we obtain
\begin{equation}\label{eq2.13}
\int\limits_a^b\E^{2\pi i \omega x}\varphi(x)\ \d x=(b-a)\E^{2\pi i \omega a}\int\limits_0^1\E^{2\pi i \omega (b-a)y}\varphi((b-a)y+a)\d y.
\end{equation}
Then, applying  (\ref{eq1.13}), Theorems \ref{Thm6} and \ref{Thm7} to the integral on the right-hand side of the last equality, we have the following results which are optimal quadrature formulas of the form (\ref{eq2.12}) in the sense of Sard in the space $W_2^{(2,1)}[a,b]$ for all real $\omega$.

For the case $\omega=0$, using (\ref{eq1.13}) in (\ref{eq2.13}), we get

\begin{theorem}\label{Thm8} Coefficients of the optimal quadrature formulas of the form
\begin{equation*}
\int\limits_a^b\varphi(x)\d x\cong \sum\limits_{\beta=0}^NC_{\beta,0}[a,b]\varphi(h\beta+a)
\end{equation*}
in the space $W_2^{(2,1)}[a,b]$ have the form
\begin{eqnarray*}
\mathring{C}_{0,0}[a,b]&=&(b-a)\left(1-\frac{\frac{1}{N}}{\E^{\frac{1}{N}}-1}-K(\lambda_1-\lambda_1^N)\right),\nonumber\\
\mathring{C}_{\beta,0}[a,b]&=&(b-a)\left(\frac{1}{N}+K\left((\E^{\frac{1}{N}}-\lambda_1)\lambda_1^{\beta}+(1-\lambda_1\E^{\frac{1}{N}})\lambda_1^{N-\beta}\right)\right),\ \\ \beta=1,2,...,N-1,\label{eq2.16}\\
\mathring{C}_{N,0}[a,b]&=&(b-a)\left(-1+\frac{\E^{\frac{1}{N}}\frac{1}{N}}{\E^{\frac{1}{N}}-1}-K(\lambda_1-\lambda_1^N)\E^{\frac{1}{N}}\right),\nonumber
\end{eqnarray*}
where
\begin{eqnarray*}
&&K=\frac{(2\E^{\frac{1}{N}}-2-\frac{1}{N}\E^{\frac{1}{N}}-\frac{1}{N})(\lambda_1-1)}{2(\E^{\frac{1}{N}}-1)^2(\lambda_1+\lambda_1^{N+1})},\quad \quad \quad \quad \quad \quad \quad\quad\quad\quad\quad\quad\quad\quad \quad \quad\quad\quad\quad\quad
\end{eqnarray*}
$$
\lambda_1=\frac{\frac{1}{N}(\E^{\frac{2}{N}}+1)-\E^{\frac{2}{N}}+1-(\E^{\frac{1}{N}}-1)\sqrt{(\frac{1}{N})^2(\E^{\frac{1}{N}}+1)^2+\frac{2}{N}(1-\E^{\frac{1}{N}})}}
{1-\E^{\frac{2}{N}}+\frac{2}{N}\E^{\frac{1}{N}}},\ \ |\lambda_1|<1.
$$
\end{theorem}

For the case $\omega\in \mathbb{R}$ and $\omega h\not\in \mathbb{Z}$, applying Theorem \ref{Thm6} to the right-hand side of (\ref{eq2.13}),
we obtain the following result.

\begin{theorem}\label{Thm9}
{For real $\omega$ with $\omega h\not\in \mathbb{Z}$, the coefficients of optimal quadrature
formulas of the form (\ref{eq2.12})  in the space
$W_2^{(2,1)}[a,b]$  are expressed by formulas}
\begin{eqnarray*}
&&\mathring{C}_{0,\omega}[a,b]=(b-a)\E^{2\pi i\omega a}\Bigg(\frac{K_{\omega,2}\ \E^{2\pi i\omega h}}{\E^{2\pi i\omega h}-1} -\frac{1}{2\pi i\omega(b-a)}+a_1\frac{\lambda_1}{\lambda_1-1}+b_1\frac{\lambda_1^{N}}{1-\lambda_1}\Bigg),\\
&&\mathring{C}_{\beta,\omega}[a,b]=(b-a)\E^{2\pi i\omega a}\Bigg(\E^{2\pi i\omega h\beta}K_{\omega,2}+a_1\lambda_1^\beta+b_1\lambda_1^{N-\beta}\Bigg),\ \quad\quad \ \beta=1,2,...,N-1,\\
&&\mathring{C}_{N,\omega}[a,b]=(b-a)\E^{2\pi i\omega a}\Bigg(K_{\omega,2}\ \Bigg(\frac{\E^{2\pi i\omega(b-a)}\E^{\frac{1}{N}}}{\E^{\frac{1}{N}}-\E^{2\pi i\omega h}}+\frac{\E^{2\pi i\omega h+1}(1-\E^{\frac{1}{N}})}{(\E^{2\pi i\omega h}-1)(\E^{\frac{1}{N}}-\E^{2\pi i\omega h})} \Bigg) \\
&&\quad\quad\quad\quad \quad \ +\frac{\E^{2\pi i\omega(b-a)}}{2\pi i\omega(b-a)-1}+a_1\frac{\lambda_1^{N}\E^{\frac{1}{N}}}{\E^{\frac{1}{N}}-\lambda_1} +b_1\frac{\lambda_1\E^{\frac{1}{N}}}{\lambda_1\E^{\frac{1}{N}}-1}\Bigg),
\end{eqnarray*}
{where }
$$
\begin{array}{l}
a_1=\frac{(\lambda_1-1)(\lambda_1-\E^{\frac{1}{N}})}{\lambda_1(1-\lambda_1^{2N})}
\Bigg(\frac{1-\lambda_1^N\E^{2\pi i\omega (b-a)}}{2\pi i\omega(b-a)(2\pi i\omega(b-a)-1)(\E^{\frac{1}{N}}-1)}+\frac{K_{\omega,2}\lambda_1^N(\E-\E^{2\pi i\omega(b-a)})\E^{2\pi i\omega h}}{(\E^{2\pi i\omega h}-1)(\E^{\frac{1}{N}}-\E^{2\pi i\omega h})}\Bigg),\\
b_1= \frac{(\lambda_1-1)(1-\lambda_1\E^{\frac{1}{N}})}{\lambda_1(1-\lambda_1^{2N})}
\Bigg(\frac{\lambda_1^N-\E^{2\pi i\omega(b-a)}}{2\pi i\omega(b-a)(2\pi i\omega(b-a)-1)(\E^{\frac{1}{N}}-1)}+\frac{K_{\omega,2}(\E-\E^{2\pi i\omega(b-a)})\E^{2\pi i\omega h}}{(\E^{2\pi i\omega h}-1)(\E^{\frac{1}{N}}-\E^{2\pi i\omega h})}\Bigg)\\
\end{array}
$$
and
$$
K_{\omega,2}=\frac{ 2(\E^{\frac{1}{N}}-\E^{2\pi i\omega h})(\E^{-2\pi i\omega h}-\E^{\frac{1}{N}})(1-\E^{2\pi i\omega h})^2}{(2\pi\omega(b-a))^2((2\pi\omega(b-a))^2+1)\left((1-\E^{\frac{2}{N}})(1-\E^{2\pi i\omega h})^2+\frac{2}{N}(\E^{\frac{1}{N}}-\E^{2\pi i\omega h})(1-\E^{2\pi i\omega h}\E^{\frac{1}{N}})\right)}.
$$
\end{theorem}

Lastly, for the case $\omega h\in \mathbb{Z}$ with $\omega\neq 0$ application of Theorem \ref{Thm7} in (\ref{eq2.13}) gives the following.

\begin{theorem}\label{Thm10} {For $\omega h\in \mathbb{Z}$ with $\omega\neq 0$, the coefficients of optimal quadrature
formulas of the form (\ref{eq2.12}) in the space
$W_2^{(2,1)}[a,b]$  are expressed by formulas}
\begin{eqnarray*}
\mathring{C}_{0,\omega}[a,b]&=&(b-a)\Bigg(\frac{2\pi i\omega(b-a) (1-\E^\frac{1}{N})-1}{2\pi i\omega(b-a) (1-2\pi i\omega(b-a))(1-\E^\frac{1}{N})}\\
&&+\frac{1}{1-\lambda_1}\Bigg(\frac{a_1\lambda_1^2}
{\E^\frac{1}{N}-\lambda_1}+\frac{b_1\lambda_1^N}
{1-\lambda_1\E^\frac{1}{N}}\Bigg)\Bigg),\\
\mathring{C}_{\beta,\omega}[a,b]&=&(b-a)\Bigg(a_1\lambda_1^\beta+b_1\lambda_1^{N-\beta}\Bigg),
\ \ \ \ \beta=1,2,...,N-1,\\
\mathring{C}_{N,\omega}[a,b]&=&(b-a)\Bigg(\frac{2\pi i\omega(b-a)(\E^\frac{1}{N}-1)-\E^\frac{1}{N}}{2\pi i\omega(b-a)(1-2\pi i\omega(b-a))(1-\E^\frac{1}{N})}\\
&&+ a_1\frac{\E^\frac{1}{N}\lambda_1^N}
{(1-\lambda_1)(\E^\frac{1}{N}-\lambda_1)}+b_1\frac{\E^\frac{1}{N}\lambda_1^2}
{(1-\lambda_1)(1-\lambda_1\E^\frac{1}{N})}\Bigg).
\end{eqnarray*}
{
where}\\
$$
a_1= \frac{(\E^\frac{1}{N}-\lambda_1)(1-\lambda_1)}{2\pi i \omega(b-a) (2\pi i \omega(b-a) -1)\lambda_1(\E^\frac{1}{N}-1)(\lambda_1^N+1)},
$$
$$
b_1=\frac{(1-\E^\frac{1}{N}\lambda_1)(1-\lambda_1)}{2\pi i \omega(b-a) (2\pi i \omega(b-a) -1)\lambda_1(\E^\frac{1}{N}-1)(\lambda_1^N+1)}.
$$
\end{theorem}

\section{Numerical results}
Here we consider some numerical results which confirm the theoretical results of the previous sections.
Here  using  the optimal quadrature formula (\ref{eq2.12}) we approximate the integrals
\begin{eqnarray}
g_{1,\omega}[-1,1]&=&\int\limits_{-1}^1\E^{2\pi \i \omega x}x\d x=\left\{
\begin{array}{ll}
\frac{2\i}{(2\pi \omega)^2}(\sin 2\pi \omega-2\pi \omega\cos 2\pi\omega),&\ \omega\neq 0,\\
0,&\ \omega=0,
\end{array}
\right. \label{eq4.1}\\
g_{2,\omega}[-1,1]&=&\int\limits_{-1}^1\E^{2\pi i \omega x}\E^{x}\d x=\left\{
\begin{array}{ll}
\frac{\E^{2\pi \i\omega+1}-\E^{-2\pi \i\omega-1}}{2\pi \i \omega+1},\ \quad  \quad \quad  \quad \quad  \quad&\omega\neq 0,\\
\frac{\E^2-1}{\E},&\ \omega=0.
\end{array}
\right. \label{eq4.2}\\
g_{3,\omega}[-1,1]&=&\int\limits_{-1}^1\E^{2\pi \i \omega x}x\E^x\d x=\left\{
\begin{array}{ll}
\frac{\E^{2\pi \i\omega+1}+\E^{-2\pi \i\omega-1}}{2\pi \i \omega+1}-\frac{\E^{2\pi \i\omega+1}-\E^{-2\pi \i\omega-1}}{(2\pi \i \omega+1)^2},&\ \omega\neq 0,\\
2\E^{-1},&\ \omega=0,
\end{array}
\right. \label{eq4.3}
\end{eqnarray}

We denote the absolute value of the error of the optimal quadrature formula (\ref{eq2.12})  by
$$
R_{\varphi,\omega}[a,b]=\left|\int\limits_a^b\E^{2\pi \i\omega x}\varphi(x)\d x-
\sum\limits_{\beta=0}^N{C}_{\beta,\omega}[a,b]\varphi(h\beta+a)\right|.
$$

Using last equality, for integrant functions $\varphi(x)=x$, $\varphi(x)=\E^x$ and $\varphi(x)=x\E^x$ when $N = 1; 10; 100$ and $\omega=1.01; 10.01; 100.01; 1000.01; 10000.01$, we get numerical results which are presented in Table 1, Table 2 and Table 3, respectively.
Numbers in parenthesis indicates the decimal exponents.

\begin{center}
\begin{table}[h]
\caption{Numerical values of $R_{x,\omega}$ for some selected values of $N$ and $\omega$  }
\centering
\begin{tabular}{l|l l l l l }
\hline
 N & $\omega=1.01$ &$\omega=10.01$& $\omega=100.01$ & $\omega=1000.01$ & $\omega=10000.01$\\ \hline
 1 &   2.436(-2)  & 2.520(-4)      & 2.527(-6)      & 2.527(-8)       & 2.528(-10)     \\
 10&   4.531(-2)  & 1.431(-5)      & 1.456(-7)      & 1.459(-9)       & 1.459(-11)     \\
$100$ &4.190(-2)  & 4.899(-5)      & 1.434(-8)      & 1.457(-10)      & 1.459(-12)      \\\hline
\end{tabular}
\end{table}
\end{center}

\begin{center}
\begin{table}[h]
\caption{Numerical values of $R_{\E^x,\omega}$ for some selected values of $N$ and $\omega$  }
\centering
\begin{tabular}{l|l l l l l }
\hline
 N & $\omega=1.01$ &$\omega=10.01$& $\omega=100.01$ & $\omega=1000.01$ & $\omega=10000.01$\\ \hline
 1 & 8.473(-2)    & 8.882(-4)   & 8.909(-6)      &  8.912(-8)       & 8.912(-10)     \\
 10& 1.791(-1)    & 6.458(-5)   & 6.584(-7)       & 6.596(-9)       & 6.597(-11)     \\
$100$ &1.706(-1)  & 1.995(-3)   & 6.622(-8)       & 6.729(-10)      & 6.741(-12)      \\\hline
\end{tabular}
\end{table}
\end{center}
\begin{center}
\begin{table}[h]
\caption{Numerical values of $R_{x\E^x,\omega}$ for some selected values $N$ and $\omega$  }
\centering
\begin{tabular}{l|l l l l l }
\hline
 N & $\omega=1.01$ &$\omega=10.01$& $\omega=100.01$ & $\omega=1000.01$ & $\omega=10000.01$\\ \hline
 1 &1.643(-1)    & 1.757(-3)   & 1.763(-5)       & 1.764(-7)     & 1.764(-9)     \\
 10& 3.688(-1)    & 1.554(-4)   & 1.587(-6)      & 1.590(-8)      & 1.590(-10)     \\
$100$ &3.584(-1) & 4.191(-3) & 1.606(-7)         & 1.633(-9)      & 1.635(-11)      \\\hline
\end{tabular}
\end{table}
\end{center}
These numerical results confirms numerical covergence of our quadrature formula.

\section{Conclusion}

In this paper we constructed optimal quadrature formulas of the form (\ref{eq1.3}) for the interval $[0,1]$ when $\omega\in \mathbb{R}$ and $\omega\neq 0$.
And we found the coefficients of the optimal quadrature formula are found. Then  using  the optimal quadrature formula (\ref{eq2.12}) we approximated the integrals.

\maketitleSecondary
\Russian

\end{document}